\documentstyle[12pt]{article}
\makeatletter
\font\twtyeuf=eufm10  \@magscale4
\font\twtymsb=msbm10 \@magscale4
\font\svteuf=eufm10  \@magscale3
\font\svtmsb=msbm10 \@magscale3
\font\frteuf=eufm10  \@magscale2
\font\frtmsb=msbm10 \@magscale2
\font\twleuf=eufm10 \@magscale1
\font\twlmsb=msbm10 \@magscale1
\font\elveuf=eufm10 \@halfmag
\font\elvmsb=msbm10 \@halfmag
\font\teneuf=eufm10
\font\tenmsb=msbm10
\font\nineuf=eufm9
\font\ninmsb=msbm9
\font\egteuf=eufm8
\font\egtmsb=msbm8
\font\seveuf=eufm7
\font\sevmsb=msbm7
\font\sixeuf=eufm6
\font\sixmsb=msbm6
\font\fiveuf=eufm5
\font\fivmsb=msbm5

\def\frak{\protect\pfrak}
\def\Bbb{\protect\pBbb}

\newfam\msbfam
\newfam\euffam

\@addfontinfo\@xxpt{
\def\pBbb{\fam\msbfam\twtymsb}\textfont\msbfam\twtymsb
     \scriptfont\msbfam\frtmsb \scriptscriptfont\msbfam\twlmsb%
\def\pfrak{\fam\euffam\twtyeuf}\textfont\euffam\twtyeuf
     \scriptfont\euffam\frteuf \scriptscriptfont\euffam\twleuf%
}
\@addfontinfo\@xviipt{
\def\pBbb{\fam\msbfam\svtmsb}\textfont\msbfam\svtmsb
     \scriptfont\msbfam\twlmsb \scriptscriptfont\msbfam\tenmsb%
\def\pfrak{\fam\euffam\svteuf}\textfont\euffam\svteuf
     \scriptfont\euffam\twleuf \scriptscriptfont\euffam\teneuf%
}
\@addfontinfo\@xivpt{
\def\pBbb{\fam\msbfam\frtmsb}\textfont\msbfam\frtmsb
     \scriptfont\msbfam\tenmsb \scriptscriptfont\msbfam\sevmsb%
\def\pfrak{\fam\euffam\frteuf}\textfont\euffam\frteuf
     \scriptfont\euffam\teneuf \scriptscriptfont\euffam\seveuf%
}
\@addfontinfo\@xiipt{
\def\pBbb{\fam\msbfam\twlmsb}\textfont\msbfam\twlmsb
     \scriptfont\msbfam\egtmsb \scriptscriptfont\msbfam\sixmsb%
\def\pfrak{\fam\euffam\twleuf}\textfont\euffam\twleuf
     \scriptfont\euffam\egteuf \scriptscriptfont\euffam\sixeuf%
}
\@addfontinfo\@xipt{
\def\pBbb{\fam\msbfam\elvmsb}\textfont\msbfam\elvmsb
     \scriptfont\msbfam\egtmsb \scriptscriptfont\msbfam\sixmsb%
\def\pfrak{\fam\euffam\elveuf}\textfont\euffam\elveuf
     \scriptfont\euffam\egteuf \scriptscriptfont\euffam\sixeuf%
}
\@addfontinfo\@xpt{
\def\pBbb{\fam\msbfam\tenmsb}\textfont\msbfam\tenmsb
     \scriptfont\msbfam\sevmsb \scriptscriptfont\msbfam\fivmsb%
\def\pfrak{\fam\euffam\teneuf}\textfont\euffam\teneuf
     \scriptfont\euffam\seveuf \scriptscriptfont\euffam\fiveuf%
}
\@addfontinfo\@ixpt{
\def\pBbb{\fam\msbfam\ninmsb}\textfont\msbfam\ninmsb
     \scriptfont\msbfam\sixmsb \scriptscriptfont\msbfam\fivmsb%
\def\pfrak{\fam\euffam\nineuf}\textfont\euffam\nineuf
     \scriptfont\euffam\sixeuf \scriptscriptfont\euffam\fiveuf%
}
\@addfontinfo\@viiipt{
\def\pBbb{\fam\msbfam\egtmsb}\textfont\msbfam\egtmsb
     \scriptfont\msbfam\sixmsb \scriptscriptfont\msbfam\fivmsb%
\def\pfrak{\fam\euffam\egteuf}\textfont\euffam\egteuf
     \scriptfont\euffam\sixeuf \scriptscriptfont\euffam\fiveuf%
}
\@addfontinfo\@viipt{
\def\pBbb{\fam\msbfam\sevmsb}\textfont\msbfam\sevmsb
     \scriptfont\msbfam\sixmsb \scriptscriptfont\msbfam\fivmsb%
\def\pfrak{\fam\euffam\seveuf}\textfont\euffam\seveuf
     \scriptfont\euffam\sixeuf \scriptscriptfont\euffam\fiveuf%
}
\@addfontinfo\@vipt{
\def\pBbb{\fam\msbfam\sixmsb}\textfont\msbfam\sixmsb
     \scriptfont\msbfam\sixmsb \scriptscriptfont\msbfam\sixmsb%
\def\pfrak{\fam\euffam\sixeuf}\textfont\euffam\sixeuf
     \scriptfont\euffam\sixeuf \scriptscriptfont\euffam\sixeuf%
}
\@addfontinfo\@vpt{
\def\pBbb{\fam\msbfam\fivmsb}\textfont\msbfam\fivmsb
     \scriptfont\msbfam\fivmsb \scriptscriptfont\msbfam\fivmsb%
\def\pfrak{\fam\euffam\fiveuf}\textfont\euffam\fiveuf
     \scriptfont\euffam\fiveuf \scriptscriptfont\euffam\fiveuf%
}

\makeatother
\newtheorem{thm}{Theorem}[section]
\newtheorem{propose}{Proposition}[section]
\newtheorem{lemma}[propose]{Lemma}
\newtheorem{cor}[propose]{Corollary}

\def\P{{\Bbb P}}   
\def\C{{\Bbb C}}     
\def\Z{{\Bbb Z}}     

\newcommand{\implies}{\mbox{$\Rightarrow$}}
\newcommand{\into}{\hookrightarrow}

\renewcommand{\iff}{\mbox{ $\Leftrightarrow$ }}
\newcommand{\ie}{{\it i.e.\/},\ }

\newcommand{\sext}{\mbox{${\cal E}xt\,$}}  
\newcommand{\limdir}[1]{{\displaystyle{\mathop{\rm
lim}_{\buildrel\longrightarrow\over{#1}}}}\,}
\newcommand{\liminv}[1]{{\displaystyle{\mathop{\rm
lim}_{\buildrel\longleftarrow
   \over{#1}}}}\,}
\newcommand{\onto}{\mbox{$\to\!\!\!\to$}}
\newcommand{\supp}{{\rm supp}\,}
\newenvironment{proof}{
                       \trivlist \item[\hskip \labelsep{\bf Proof}:]
                      }{
                        \hfill$\Box$\endtrivlist
                      }
\newenvironment{rmk}{
                        \trivlist \item[\hskip \labelsep{\bf Remark}:]
                      }{
                        \endtrivlist}

\newcommand{\rest}{\mbox{${\cal O}(D)_{\mid D}$}}
\newcommand{\resta}[1]{\mbox{${\cal O}({#1D})_{\mid D}$}}
\newcommand{\restb}[2]{\mbox{${\cal O}({#1D})_{\mid{#2D}}$}}
\newcommand{\lrw}{\mbox{$\longrightarrow$}}
\newcommand{\cO}{\mbox{$\cal O$}}
\newcommand{\ohd}[1]{\mbox{${\rm H}^0(\cO_{#1D})$}}
\newcommand{\fom}[1]{\mbox{$\widehat{\Omega}^{#1}_{\frak X}$}}
\newcommand{\fo}{\mbox{$\widehat{\cal O}_{\frak X}$}}

\title {Affine like Surfaces }
\author{N. Mohan Kumar\\
                School of Mathematics\\
                Tata Institute of Fundamental Research\\
                Homi Bhabha Road, Bombay 400 005}

\begin{document}
\maketitle
\begin{center}
({\it Dedicated to the memory of my parents})
\end{center}
\section{Introduction}
If $X$ is a variety, then $X$ is affine if and only if
${\rm H}^i(X, {\cal F})=0$ for
all coherent sheaves $\cal F$ and for all positive $i$.
This paper deals with the following natural question:

\noindent{\it Question}: Classify all smooth varieties $X$ (over $\C$)
with ${\rm H}^i(X, \Omega^j_X)=0$ for all $j$ and for all
positive $i$.

Of course if $\dim X=1$ such an $X$ is affine. Here we deal with
the case of surfaces and completely classify them. This question
was raised by T.~Peternell \cite{Peternell}. Our theorem is
as follows:

\begin{thm} Let $Y$ be a smooth surface with ${\rm H}^i(Y,
\Omega^j)=0$ for all $j$ and $i>0$. Then $Y$ is one of the
following.
\begin{enumerate}
\item $Y$ is affine.
\item Let $C$ be an elliptic curve and $E$ the unique non-split
extension of $\cal O_C$ by itself. Let $X=\P(E)$ and $D$ the
canonical section. $Y=X-D$.
\item Let $X$ be a projective rational surface with an effective
divisor $D=-K$ with $D^2=0$, \rest is non-torsion
and the dual graph of $D$ be $\tilde{D}_8$ or $\tilde{E}_8$. $Y=X-D$.
\end{enumerate}
\end{thm}

\noindent{\it Remark}. I do not know whether the surfaces in the
last case of the theorem are Stein. It is well known that in the
second case they are Stein.

\noindent{\it Acknowledgements}. The paper of Peternell was
brought to my attention by R.~R.~Simha.
I thank him. I have benefitted by many
discussions I had with V.~Srinivas and K.~Paranjape. In
particular, when I was unable to settle some of the cases, it
was Srinivas who suggested that I use formal de Rham cohomology
and showed me how. I thank both of them.

\section{Some preliminary lemmas}
We assume that $Y$ is a smooth surface, satisfying the
hypothesis of the Theorem. Most of the lemmas in this section are
standard and many of the statements and proofs can be found in
the literaure. We include them here only for completeness.
See also Peternell's paper \cite{Peternell}.

\begin{lemma}\label{one}
 $Y$ contains no complete curves.
\end{lemma}

\begin{proof}
If $C\subset Y$ is a complete irreducible curve, then one has an
exact sequence,
$$0\lrw A\lrw \Omega^1_Y\lrw \Omega^1_C\lrw 0,$$
where $A$ is coherent. Since ${\rm H}^2(\Omega^2_Y)=0$, $Y$ is not
projective. If $X$ is a smooth projective completion of $Y$ and
if $Z=X-Y$, for any coherent sheaf $\cal F$ on $X$ we have an
exact sequence,
$${\rm H}^2_Z({\cal F})\stackrel{\alpha}{\to} {\rm H}^2(X, {\cal F})\to {\rm
H}^2(Y, {\cal F})\to 0.$$
By  formal duality (see \cite{Rhart}
Chapter~\uppercase\expandafter{\romannumeral 3}, Theorem 3.3)
$${\rm H}^2_Z({\cal F})^*\cong {\rm H}^0(\widehat{X},
\widehat{{\cal G}})$$
where ${\cal G}={\rm Hom} ({\cal F}, \omega_X)$,  $\omega_X$
is the dualising sheaf (which in our case is just the
canonical bundle) and $`\widehat{}'$ denotes completion along $Z$.
Since $Z\neq\emptyset$, for any locally free sheaf $\cal G$
the canonical map ${\rm H}^0(X, {\cal G})\to {\rm
H}^0(\widehat{X},\widehat{{\cal G}})$ is injective
and hence $\alpha$ above  which is just the dual of this map is
surjective for locally free sheaves and thus ${\rm H}^2(Y, {\cal
F})=0$ for locally free sheaves. Since any coherent sheaf is the
quotient of a locally free sheaf we see that ${\rm H}^2(Y, {\cal
F})=0$ for all coherent sheaves and in particular ${\rm H}^2(Y, A)=0$.
Since ${\rm H}^1(\Omega^1_Y)=0$ by hypothesis, we get
${\rm H}^1(\Omega^1_C)=0$ which is absurd.
\end{proof}

\begin{lemma}\label{two}
Let $X$ be a smooth projective completion of $Y$. Then $X-Y$ is a
union of divisors.
\end{lemma}

\begin{proof}
First since $Y$ contains no complete curves, $X-Y$ cannot be of
dimension zero. We want to show that $X-Y$ has no isolated points.
If $P$ is such an isolated point, let $Y'=Y\cup \{P\}$. So $Y'$ is not
complete and $Y=Y'-\{P\}$. Writing the local cohomology exact sequence
for $\cal O$ one has,
$$0={\rm H}^1(Y, {\cal O}_Y)\to {\rm H}^2_{\{P\}}({\cal O})\to
{\rm H}^2(Y', {\cal O}_{Y'})=0$$
The first zero by hypothesis and the last zero since $Y'$ is not complete.
But the middle term is infinite dimensional, since
$${\rm H}^2_{\{P\}}({\cal O})=\limdir{n} {\rm Ext}^2({\cal
O}_{nP}, {\cal O}_{Y'})
=\limdir{n} {\rm H}^0({\cal E}xt^2({\cal O}_{nP}, {\cal O}_{Y'}))$$
where ${\cal O}_{nP}={\cal O}/{\frak M}^n$, $\frak M$ the ideal sheaf
defining $P\in Y'$. This leads us to a contradiction.
\end{proof}

We may assume that $X-Y=\cup_{i=1}^n E_i$ and no $E_i$ is
exceptional of the first kind. From now on we will also assume
that $Y$ is not affine.

\begin{lemma}\label{three}
Let $D$ be any effective divisor on $X$ with $\supp D=\cup E_i$.
If $V$ is a vector bundle on $X$, then
$${\rm H}^1(Y, V_{\mid Y})=0 \iff \limdir{n}{\rm H}^1(X, V(nD))=0.$$
\end{lemma}

\begin{proof}
Writing the local cohomology sequence we see that ${\rm H}^1(Y,
V_{\mid Y})=0$ is equivalent to,

$$\limdir{n} {\rm H}^1(\sext^1(\cO_{nD}, V)= {\rm H}^2(X, V)
\ {\rm and}\ \limdir{n} {\rm H}^0(\sext^1(\cO_{nD}, V))\onto
{\rm H}^1(X, V).$$

We have $\sext^1(\cO_{nD}, V)=V(nD)_{\mid nD}.$
We have a commutative diagram,
\[
\begin{array}{rcccccccl}
0&\to& V&\to&V(nD)&\to& V(nD)_{\mid nD}&\to& 0\\
&&\mid\mid&&\downarrow&& \downarrow&&\\
0&\to& V&\to& V((n+l)D)&\to&V((n+l)D)_{\mid (n+l)D}&\to& 0
\end{array}
\]

This gives the following exact sequence,
$$\limdir{n}{\rm H}^0(V(nD)_{\mid nD})\stackrel{\alpha}{\to}
{\rm H}^1(X, V)\to
\limdir{n}{\rm H}^1(V(nD))\to$$
$$\limdir{n}{\rm H}^1(V(nD)_{\mid nD})\stackrel{\beta}{\to} {\rm H}^2(X,
V)\to \limdir{n}{\rm H}^2(V(nD))=0$$

{}From this  it follows that
$\limdir{n}{\rm H}^1(V(nD))$ is zero if and only $\alpha$ is
surjective and $\beta$ is an isomorphism which in turn is
equivalent to ${\rm H}^1(Y, V_{\mid Y})=0$.
\end{proof}

\begin{lemma}\label{four}
$X-Y=\cup_{i=1}^n E_i$ is connected.
\end{lemma}

\begin{proof}
If $C\subset X$ is any effective divisor,
${\rm H}^2_C(K)\onto {\rm H}^2(X, K)=\C$
where $K$ is the canonical bundle. If $X-Y=D=D_1+D_2$ and $D_1\cap
D_2=\emptyset$ then ${\rm H}^2_D(K)={\rm H}^2(K)=\C$ since ${\rm
H}^i(Y, K_Y)=0$ for $i=1,2$; but
${\rm H}^2_D(K)={\rm H}^2_{D_1}(K)\oplus {\rm H}^2_{D_2}(K)$,
both summands of which are at least one
dimensional. This is a contradiction.
\end{proof}

Next we analyse the intersection form of the $E_i$'s.

\begin{lemma}\label{five}
If $D$ is any divisor with $\supp D\subset\cup E_i$ then
$D^2\leq 0$.
\end{lemma}

\begin{proof}
If not there exists a $G=\sum a_iE_i$ with $G^2>0$. Writing
$G=G_1-G_2$, with $G_i$'s effective,we see that
$0<G^2=G_1^2-2G_1G_2+G_2^2$. Since $G_1G_2\geq 0$ we see that
$G_1^2>0$ or $G_2^2>0$. So we may assume that $G$ is effective.
Let $G=G_1+G_2$ be the Zariski decomposition of $G$. Then $G_1$
is arithmetically effective and for any curve $E\subset G_2$,
$G_1\cdot E=0$. So $0<G^2=G_1^2+G_2^2$. The intersection form on
$\supp G_2$ is negative definite and so $G_1^2>0$. Thus we may
assume that $G$ is effective and $G\cdot E\geq 0$ for all curves
$E\subset\supp G$ (actually we may assume that $G$ is
arithmetically effective).

Now let $E_1, \ldots , E_k$ be the divisors in $\supp G$ such
that $G\cdot E_i=0$. Write $G=A+B$ with $B$ consisting of all
the $E_i$'s. We will show that there exists an effective
divisor $G'$ which has the property that $G'\cdot E>0$ for every
$E\subset\supp G'$. The proof is by induction on $k$, the number of
components of $B$. For any $E_i$ since $G\cdot E_i=0$ we see
that $B\cdot E_i\leq 0$. If $B\cdot E_i=0$ for all $i$, then
$\supp A\cap\supp B=\emptyset$ and hence $A$ will do the job. So
assume that there exists at least one $E_i$ with $E_i\cdot B<0$,
say $E_1$.  Consider $G'=mG-B$. For large $m$, $G'$ is effective
and for all $E\subset\supp A$, $G'\cdot E>0$. Also $G'\cdot
E_i=-B\cdot E_i\geq 0$. But $E_1\cdot G'=-E_1\cdot B>0$. This
gives the induction step.

Next we show that $\supp G$ can be taken to be $\cup_{i=1}^n
E_i$. If there is an $E$ outside the support of $G$, since
$\cup_{i=1}^n E_i$ is connected there is such an $E$ with
$G\cdot E>0$. Now consider $G'=mG+E$ for large $m$. Then
$G'^2>0$, $G'\cdot E'>0$ for all $E'\subset\supp G$ and $G'\cdot
E>0$. Thus by an easy induction we are done. If $C$ is any curve
we get that $G\cdot C>0$ since $Y$ contains no complete curves.
Thus $G$ is ample, by the Nakai-Moisezhon criterion. Hence $Y$
is affine, a contradiction.
\end{proof}

\begin{lemma}\label{six}
The intersection form is not negative definite. There exists a
unique (effective) divisor $D=\sum a_iE_i$ with $gcd(a_i)=1$
generating the kernel of this intersection form.
\end{lemma}

\begin{proof}
If the intersection form is negative definite, then by
corollary 2.11, Chapter \uppercase\expandafter{\romannumeral 1},
\S 2 \cite{Barth}  we get an effective
divisor $D=\sum a_iE_i$ with
$a_i>0$ and $D\cdot E_i<0$ for all $i$. It is easy to see from
this that
$$h^0(\cO_D(nD))=0 \ {\rm for}\ n>0 \ {\rm and}\
h^1(\cO_D(nD))\to\infty\ {\rm as}\ n\to\infty.$$
Thus from the
exact sequence
$$0\to\cO((n-1)D)\to\cO(nD)\to\cO_D(nD)\to 0,$$
one gets ${\rm H}^1(\cO((n-1)D))\into {\rm H}^1(\cO(nD))$ for all
$n>0$ and ${\rm H}^1(\cO(nD))\neq 0$ for large $n$. Thus $\limdir{n}
{\rm H}^1(nD)\neq 0$ contradicting our hypothesis that ${\rm
H}^1(Y, {\cal O}_Y)=0$ by lemma [\ref{three}].

Thus the intersection form is not negative definite. Then by
Chapter \uppercase\expandafter{\romannumeral 5}, \S 3.5
\cite{Bourbaki}  we get the result.
\end{proof}

{}From now on we fix such a divisor $D=\sum a_iE_i$, as the
generator of the kernel of the intersection form.

\begin{cor}\label{sixh}
Let $G$ be any divisor with $\supp G\subset\supp D$ and assume
that $G\cdot E=0$ for all $E\subset\supp G$. Then $G=nD$ for
some $n$.
\end{cor}

\begin{proof}
If $G\cdot E=0$ for all $E\subset\supp D$ we are done by the lemma.
So assume that there exists an $E\subset\supp D$ such that $G\cdot
E\neq 0$. Then for any $n$, we have
$$(nG+E)^2=2n(G\cdot  E)+E^2$$ and this can be made positive by
choosing $n$ sufficiently large or small according to the sign
of $G\cdot E$. This contradicts lemma [\ref{five}].
\end{proof}

\begin{lemma}\label{seven}
${\rm H}^0({\cal O}_X(mD))=\C$ for all $m$.
\end{lemma}

\begin{proof}
If $h^0(mD)>1$ for some $m$, write the Zariski decomposition
(see \cite{Zariski})
$mD=D_1+D_2$. If $E\subset\supp D_2$ then $D_1\cdot E=0$. If
$E\subset\supp D-\supp D_2$, then since $E\cdot D=0$ and $E\cdot
D_2\geq 0$ we see that $D_1\cdot E\leq 0$. But $D_1$ is
arithmetically effective and hence $D_1\cdot E=0$. Thus
$D_1\cdot E=0$ for every $E\subset\supp D$. Then $D_1=qD$ for
some rational number $q$. Then $D_2=0$ and hence for large $m$,
$D$ has no base components and since $D^2=0$ it has no base
points by a theorem of Zariski \cite{Zariski}. If $C\in
\mid\!mD\!\mid$, a general member, then $C$ does not meet $D$
which is not possible by lemma \ref{one}.
\end{proof}

\begin{lemma}\label{eight}
${\rm H}^0(\cO_D)=\C$.
\end{lemma}

\begin{proof}
If $D$ has only one component, this is trivial. If
$h^0(\cO_D)>1$, we can choose a maximal divisor $D_1\subset D$
such that if we write $D=D_1+D_2$, then
$h^0(\cO_{D_2}(-D_1))\neq 0$. Let $E\subset\supp D_2$. By maximality we
get $h^0(\cO_{D_2-E}(-D_1-E))=0$ and hence $h^0(\cO_E(-D_1))\neq
0$. Therefore $D_1\cdot E\leq 0$. This implies $D_2\cdot E\geq
0$ since $D\cdot E=0$. But $D_2^2\leq 0$ by negative
semi-definiteness. So $D_2\cdot E=0$ for all $E\subset\supp D_2$
and thus $D_2=nD$, $n\in\Z$ by corollary [\ref{sixh}]. This is a contradiction.
\end{proof}

\begin{lemma}\label{nine}
Let $L\in {\rm Pic}^0 D$, \ie $L$ is a line bundle on $D$ which has
degree zero when restricted to each irreducible component of $D$.
Then ${\rm H}^0(L)\neq 0$ if and only if $L\cong\cO_D$.
\end{lemma}

\begin{proof}
Let $0\neq s\in {\rm H}^0(L)$. If $s$ did not vanish on any component
of $D$ we would be done. So choose $D_1$ maximal such that
$D=D_1+D_2$ and $s=0$ in $L_{\mid D_1}$. So $s\in
{\rm H}^0(L\otimes\cO_{D_2}(-D_1))$. As before for any $E\subset\supp D_2$, we
see that ${\rm H}^0(L\otimes\cO_E(-D_1))\neq 0$ by maximality of
$D_1$. Thus $D_1\cdot E\leq 0$ which in turn implies that
$D_2\cdot E\geq 0$. Again by lemma [\ref{five}] we see that
$D_2\cdot E=0$ for all such $E$.  So again by corollary
[\ref{sixh}], $D_2=nD$ for some $n$.
\end{proof}

\begin{lemma}\label{can}
If $K$ denotes the canonical bundle of $X$, then $K\cdot E\geq 0$
for all $E\subset\supp D$.
\end{lemma}

\begin{proof}
If $\supp D$ contains more than one curve, then since $D\cdot
E=0$ for all curves in $\supp D$ and $\supp D$ is connected we
see that $E^2<0$ for all these curves. If $K\cdot E<0$ for one
of these, then it would be exceptional of the first kind which
we have assumed to be not the case. So if the lemma is false,
then $D$ must be irreducible. Since $D^2=0$ we see that
$D\cong\P^1$ and $D\cdot K=-2$. By the Riemann-Roch formula one
sees that $h^0(nD)\to\infty$ as $n\to\infty$. So choose the
smallest integer $n$ such that $h^0(nD)>1$. Then we get a $G\sim
nD$ and by our choice of $n$, $G\cap D=\emptyset$. This
contradicts lemma [\ref{one}].
\end{proof}

\section{Torsion case}

{}From the previous lemma, we have two possibilities. Either
$L=\cO_D(D)$ is torsion or non-torsion. We will separately
analyse the two cases. In this section we will look at the case
when it is torsion. Let $p=ord\ L$.

\begin{lemma}\label{ten}
Let $m\geq 1$ be any integer and assume that for all $l$ with
$1\leq l<m$, $\cO_{lD}(D)$ is $p$-torsion. Then $h^0({\cal O}_{mD})=a$ where
$m=ap-b$ with $0\leq b<p$.
\end{lemma}

\begin{proof}
If $m\leq p$, we must show that $h^0({\cal O}_{mD})=1$.
For $m=1$ this is just lemma [\ref{eight}].
We have
$$0\to L^{1-m}\to\cO_{mD}\to\cO_{(m-1)D}\to 0.$$
If $m>1$, then $L^{1-m}$ is not trivial since $p>m-1\geq 1$. So
$$1\leq h^0(\cO_{mD})\leq h^0(\cO_{(m-1)D})=1$$ by induction. So
we are done.

Now assume that $m>p$. We have
$$0\to\cO_{(m-p)D}(-pD)=\cO_{(m-p)D}\to\cO_{mD}\to\cO_{pD}\to 0.$$
Since ${\rm H}^0(\cO_{pD})=\C$, ${\rm H}^0(\cO_{mD})\onto {\rm
H}^0(\cO_{pD})$. So
$$h^0(\cO_{mD})=h^0(\cO_{(m-p)D})+1=a$$
by induction.
\end{proof}

\begin{lemma}\label{eleven}
Let $m$ be as in the previous lemma. Then for any $s\leq t\leq
m$ one has ${\rm H}^0(\cO_{tD})\onto {\rm H}^0(\cO_{sD})$.
\end{lemma}

\begin{proof}
If $s\leq p$, then $\ohd{s}=\C$ and `1' in $\ohd{t}$ maps to
`1'in $\ohd{s}$. So assume that $s>p$. Consider,
\[
\begin{array}{rcccccccccl}
0&\to&\cO_{(t-p)D}&=&\cO_{(t-p)D}(-pD)&\to&\cO_{tD}&\to&
\cO_{pD}&\to&0\\
&&\downarrow&&&&\downarrow&&\mid\mid&&\\
0&\to&\cO_{(s-p)D}&=&\cO_{(s-p)D}(-pD)&\to&\cO_{sD}&\to
&\cO_{pD}& \to &0
\end{array}
\]
By induction $\ohd{(t-p)}\onto \ohd{(s-p)}$. Since $\ohd{p}=\C$ ,
the map $\ohd{t}\to\ohd{p}$ is surjective. This implies the surjectivity of
$\ohd{t}\to \ohd{s}$.
\end{proof}

\begin{lemma}\label{twelve}
Let $m$ and $p$ be as above. Then $\cO_{mD}(D)$ is either
$p$-torsion or non-torsion.
\end{lemma}

\begin{proof}
$\ohd{m}\onto\ohd{(m-1)}$ by the previous lemma. Thus
${\rm H}^0(\cO_{mD}^{*}) \to {\rm H}^0(\cO_{(m-1)D}^{*})$ is
also surjective.
 Thus we get
 $$0\to {\rm H}^1(\cO_D((1-m)D))\to {\rm Pic}\ mD\to
{\rm Pic}\ (m-1)D.$$
Since ${\rm H}^1(\cO_D((1-m)D)) $ is a vector space over a field
of characteristic zero, the proof is clear.
\end{proof}

\begin{lemma}\label{thirteen}
$\cO_{mD}(D)$ is not torsion for some $m$.
\end{lemma}

\begin{proof}
If not then by the previous lemma it is $p$-torsion for all $m$.
Consider for a fixed $a$,
$$0\to\cO(aD)\to \cO(rpD)\to\cO_{(rp-a)D}(rpD)=\cO_{(rp-a)}\to
0$$ with $r>>0$. Since ${\rm H}^0(mD)=\C$ for all $m$ we see that
$\ohd{(rp-a)}\into {\rm H}^1(\cO(aD))$ for all $r>>0$. By
lemma [\ref{six}] $h^0(\cO_{(rp-a)D})\to\infty$ as $r\to\infty$. But
$h^1(aD)$ is constant independent of $r$. This is impossible.
\end{proof}

Let $k$ be the integer such that $\cO_{lD}(D)$ is $p$-torsion
for $l<k$ and $\cO_{kD}(D)$ is not torsion. By hypothesis $k\geq
2$.

\begin{lemma}\label{fourteen}
If $l\geq k$ and $m\geq l-k+1$, then
$h^0(\cO_{ld}(mD))=h^0(\cO_{(k-1)D}(m-l+k-1)D))$.
\end{lemma}

\begin{proof}
We have
$$0\to\cO_{(l-1)D}((m-1)D)\to\cO_{lD}(mD)\to\cO_D(mD)\to 0$$
If ${\rm H}^0(\cO_{(l-1)D}((m-1)D))={\rm H}^0(\cO_{lD}(mD))$, we would be
done by induction. If not ${\rm H}^0(\cO_{lD}(mD))\to{\rm H}^0(\cO_D(mD)$ is
non-zero. Then $\cO_D(mD)\cong\cO_D$ and since ${\rm H}^0(\cO_D)=\C$
we see that $\cO_{lD}(mD)\cong\cO_{lD}$. This is contrary to our
hypothesis that $l\geq k$.
\end{proof}

\begin{lemma}\label{fifteen}
$K\cdot E=0$ for $E\subset\supp D$.
\end{lemma}

\begin{proof}
By lemma [\ref{can}] we already know that $K\cdot E\geq 0$ for
all $E\subset\supp D$.
Thus to prove the lemma, it suffices to show that $K\cdot D=0$.
Consider
$$0\to\cO(mD)\to\cO((m+l)D)\to\cO_{lD}((m+l)D)\to 0.$$
and take $l$ sufficiently large and $l\geq k$. Then we get,
$h^0(\cO_{lD}((m+l)D))= h^1(\cO(mD))$, by lemma [\ref{three}].
$h^0(\cO_{lD}((m+l)D))=h^0(\cO_{(k-1)D}((m-k+1)D))$ by the
previous lemma. Since $\cO_{(k-1)D}(D)$ is $p$-torsion,
$h^0(\cO_{(k-1)D}((m-k+1)D))$ can take only finitely many values
as $m$ varies and hence
$h^1(mD)$ is bounded. Now the Riemann-Roch theorem implies that $K\cdot
D=0$.
\end{proof}

\begin{lemma}\label{sixteen}
$p_g=0$, $q=1$ and $D$ is smooth and irreducible.
\end{lemma}

\begin{proof}
For large $s$, ${\rm H}^1(\cO_X)\to {\rm H}^1(\cO(sD))$ is zero
by lemma [\ref{three}]. So we get ${\rm H}^0(\restb{s}{s})={\rm
H}^1(\cO_X)$ by lemma [\ref{seven}]. If further $s\geq k$, then
by lemma [\ref{fourteen}], we get $q=h^0(\restb{(k-1)}{(k-1)})$.
$$0\to\cO(mD)\to\cO((m+l)D)\to\restb{(m+l)}{l}\to 0$$ gives, for
sufficiently large $l\geq k$,
$$h^1(mD)=h^0(\restb{(m+l)}{l})=h^0(\restb{(m+k-1)}{(k-1)})$$
by lemma [\ref{fourteen}]. If $m=ap$, we get that
$h^1(apD)=h^0(\restb{(k-1)}{(k-1)}) = q$, for all $a>>0$. By the
Riemann-Roch theorem we have, $1-h^1(apD)=1-q+p_g$ and thus $p_g=0$.
Thus for all large $l$, $h^1(lD)=q$. If $q=0$, then let $l=ap$,
$a>>0$. We have from the exact sequence,
$$0\to\cO((ap-1)D)\to\cO(apD)\to\cO_D\to 0$$ since ${\rm
H}^1((ap-1)D)=0$ for large $a$, $h^0(apD)\geq 2$, contradicting
lemma [\ref{six}]. So $q\geq 1$.

Consider $X\stackrel{\pi}{\lrw}\ {\rm Alb}\ X$. Since $p_g=0$,
$\pi(X)$ is a curve. Since $Y$ contains no complete curves by
lemma [\ref{one}], $\pi(D)$ cannot be a point. If the genus of
$\pi(X)\geq 2$ then $\pi(D)$ is a point since any component of
$D$ has arithmetic genus at most one. So $q=1$ and
$D\lrw\ {\rm Alb}\ X$ is surjective. If $D$ is not irreducible
then each component $E$ of $D$ has self-intersection negative
since $D\cdot E=0$. Since $E\cdot K=0$, $E\cong\P^1$ and then
$\pi(D)$ is a point. So $D$ is irreducible. If $D$ is not
smooth, since its arithmetic genus is one, its geometric genus
must be zero and hence again $\pi(D)$ will be a point. So $D$
must be  a smooth elliptic curve.
\end{proof}

\begin{lemma}\label{seventeen}
$X$ is ruled.
\end{lemma}

\begin{proof}
If not $h^0(mK)\neq 0$ for some $m>0$. Write $mK=\sum n_iF_i$.
If $F_i\neq D$ then since $D\cdot K=0$, $F_i\cap D=\emptyset$
and hence $F_i\subset Y$ which is not possible. So $mK=nD$ for
some integer  $n$ and if the surface is not ruled, then
$n\geq 0$. Hence $K^2=0$ and ${\rm H}^0({\cal O}(-aK-bD))=0$ for
all $a,b>0$. For any positive
integer $l$ we have an exact sequence,
$$0\to\cO(lK+(l-1)D)\to\cO(lK+lD)\to\cO_D\to 0$$
If $l\geq 2$ then by duality, ${\rm H}^2(lK+(l-1)D)={\rm H}^2(lK+lD)=0$
and hence
${\rm H}^1(lK+lD)\onto {\rm H}^1(\cO_D)=\C$. By the Riemann--Roch
theorem one
sees that $h^0(lK+lD)\geq 1$ for all $l\geq 2$. For any such $l$ let
us write $lK+lD=D_l=a_lD+G$ with $D$ not in the support of $G$.
Since the complement of $D$ contains no complete
curves while $D\cdot G=0$, $G$ must be zero. Since
$mK=nD$ with $m>0$ and $n\geq 0$ one sees that $a_{l+1}>a_l$. Then
$$K+D=D_{l+1}-D_l=(a_{l+1}-a_l)D$$
and hence $K$ is effective. This is a contradiction.
\end{proof}

\begin{lemma}
$X$ is  minimal.
\end{lemma}

\begin{proof}

We have an exact sequence,
$$0\to\cO(2K+D)\to\cO(2K+2D)\to\cO_D\to 0$$
If ${\rm H}^2(2K+D)\neq 0$
then by duality ${\rm H}^0(-K-D)\neq 0$.  As before we see that
$-K-D=aD$ for some non-negative integer $a$.  Then $K^2=0$ which
implies, (since the minimal model also has $K^2\leq 0$ by Noether's
formula) that $X$ is minimal. So assume that ${\rm H}^2(2K+D)=0$ which
in turn implies that ${\rm H}^2(2K+2D)=0$ and ${\rm H}^1(2K+2D)\onto
{\rm H}^1(\cO_D)=\C$. Thus ${\rm H}^1(2K+2D)\neq 0$.  Now by
Riemann--Roch theorem one sees that ${\rm H}^0(2K+2D)\neq 0$. As before
$2K+2D=aD$ for some non-negative integer $a$. Thus again $K^2=0$ and
hence $X$ is minimal.
\end{proof}

Now let $X=\P_E(V)$, $V$ a rank two vector bundle on $E$.

\begin{lemma}
$V$ is indecomposable.
\end{lemma}

\begin{proof}
If $V$ is decomposable, we may assume that $V=\cO\oplus L$ with
$\deg L\leq 0$. Let $G$ be the section corresponding to
$\cO\into V$. Then $G_{\mid G}\cong L$. Let $F$ be a fiber.
Write $D\equiv aF+bG$. Then $a\geq 0$ and $a>0$ if $G\neq D$.
But $G\neq D$ since there is a disjoint section. So $a>0$.
$0<G\cdot D=a+b\deg L$. Also $b>0$. $0=D^2=2ab+b^2\deg L$. So
$2a+b\deg L=0$. But $2a+2b\deg L>0\implies b\deg L>0$, a
contradiction.
\end{proof}

Finally we see that $X=\P_E(V)$ with a non-split exact sequence
$$0\to\cO_E\to V\to L\to 0.$$
$\deg L=0$ or 1. It was shown by A. Neeman in \cite{Neeman} that $\deg L=1$
cannot occur. So $\deg L=0$. But the sequence does not split
implies that $L=\cO_E$. Let $G$ be the canonical section and
write $D\equiv aF+bG$ where $F$ is a fiber as before. If $D\neq
G$, then $0<G\cdot D=a$ and $0<F\cdot D=b$. But then $0=D^2=2ab$
a contradiction. So $D$ must be the canonical section.

It is easy to check that in the above situation, ${\rm H}^i(Y,
\Omega^j)=0$ for $i>0$ and $j\geq 0$.

\section{Non-torsion case}

In this section we will analyse the case when $D_{\mid D}$ is
non-torsion. So from now on let us assume this.

\begin{lemma}
For any $E\subset\supp D$, $K\cdot E=0,\ D=-K$ and $X$ is
rational.
\end{lemma}

\begin{proof}
{}From the exact sequence,
$$0\to\cO(mD)\to\cO((m+1)D)\to\resta{(m+1)}\to 0,$$
since ${\rm H}^0(\resta{(m+1)})=0$ for every $m\geq 0$ by lemma [\ref{nine}],
we see that
${\rm H}^1(mD)\into {\rm H}^1((m+1ª)D)$. But since the direct
limit of these inclusions is zero we see that ${\rm H}^1(mD)=0$.
Now by the Riemann--Roch theorem we see that $K\cdot D=0$.
Since $D$ is connected and contains no exceptional curve of the
first kind this implies that $K\cdot E=0$ for all $E\subset\supp
D$. Also ${\rm H}^1(\cO)=0$. Consider,
$$0\to K\to K+D\to K_D\to 0.$$
Since $q=0$, we get ${\rm H}^0(K+D)\onto {\rm H}^0(K_D)$. Since
$\chi(\cO_D)=0$ and $h^0(\cO_D)=1$, we see that $h^0(K_D)=1$. So
by lemma [\ref{eight}], $K_D\cong\cO_D$ and the section `1'
lifts to ${\rm H}^0(K+D)$. So $K+D=\sum b_iF_i$ where the $F_i$'s are
disjoint from $D$. Since there are no such curves $F_i$ by lemma
[\ref{one}], we see that $K=-D$. So $p_n(X)=h^0(-nD)=0$ and thus
$X$ is rational.
\end{proof}

\begin{lemma}
$D=\tilde{D}_8$ or $\tilde{E}_8$.
\end{lemma}

\begin{proof}
We have the de Rham complex,
$$0\to\cO_Y\to\Omega^1_Y\to\cO_Y\to 0$$
By a theorem of Grothendieck, \cite{Grothendieck}, we see that
since ${\rm H}^0(\cO_Y)=\C$, ${\rm H}^2(Y, \C)$ is at most
one-dimensional. We have a long exact sequence,
$$0={\rm H}^1(X, \C)\to {\rm H}^1(D, \C)\to {\rm H}^2_c(Y)\to
{\rm H}^2(X, \C)\to {\rm H}^2(D,\C)\to$$
$${\rm H}^3_c(Y)\to {\rm H}^3(X, \C)=0.$$
Since the components of $D$ are linearly independent in ${\rm Pic}\
X$, (otherwise ${\rm H}^0(Y, \cO_Y)\neq\C$), we see that ${\rm
H}^2(X, \C)\onto {\rm H}^2(D, \C)$. ${\rm H}^2(X, \C)=\C^{10}$ since $K^2=0$.
${\rm H}^3_c(Y)=0\implies {\rm H}_1(Y)=0\implies {\rm
H}^1(Y)=0$. If ${\rm H}^2(X)\to {\rm H}^2(D)$ is an isomorphism
then $D$ will support a divisor with positive self-intersection,
which we know is
not the case. Thus ${\rm H}^2_c(Y)\neq 0$.  But ${\rm
H}^2_c(Y)={\rm H}_2(Y)={\rm H}^2(Y)$ has dimension atmost one.
So ${\rm H}^2(Y)=\C$ and ${\rm H}^0(\Omega^1_Y)=0$. Also ${\rm
H}^1(D, \C)=0$ and $D$ has exactly nine components. Then by
\cite{Demazure} we see that $D$ is either $\tilde{D}_8$ or
$\tilde{E}_8$.
\end{proof}

\begin{thm}
Let $X$ be a rational surface with an effective divisor
$D\subset X$, $D=\tilde{D}_8$ or $\tilde{E}_8$, $D=-K$, $D^2=0$
and $D_{\mid D}$ is non-torsion. Let $Y=X-D$. Then ${\rm H}^i(Y,
\Omega^j_Y)=0$ for all $j$ and all $i>0$.
\end{thm}

\begin{proof}
Since $Y$ is non-complete, ${\rm H}^2(Y, \Omega^j)=0$ for all $j$.
We will show that ${\rm H}^1(X, nD)=0 \ \forall n\geq 0$ and
$\limdir{n}{\rm H}^1(X, \Omega^1(nD))=0$. The statement for $K$
follows by duality since $D=-K$. Since $D_{\mid D}$ is non-torsion one sees
immediately that ${\rm H}^0(nD)=\C$ for all $n\geq 0$. The Riemann--Roch
theorem then gives that ${\rm H}^1(nD)=0$.

Let $\frak X$ denote the formal completion of $X$ along $D$. We
have the formal de Rham complex,
$$0\to\fo\to\fom{1}\to\fom{2}\to 0$$
The hypercohomology of this complex computes the singular cohomology of
$D$ (see Chapter 4, Theorem 1.1, \cite{Hartshorne}). We have ${\rm
H}^1(D, \C)=0$ and
${\rm H}^2(D, \C)=\C^9$. We have a spectral sequence
$$E_1^{p,q}={\rm H}^q(\fom{p})\implies {\rm H}^{p+q}(D).$$
First I claim that $E_1^{p,q}=0$ if $p<0$ or $q<0$ or $p\geq 2$ or
$q\geq 2$. The first two cases are obvious.
$${\rm H}^q(\fom{p})= \liminv{n}{\rm H}^q(\Omega^p_{\mid nD})$$
and if $q\geq 2$ each term on the right is zero, since $\dim D=1$.
So the only case left is when $p\geq 2$. Of course it is trivial when
$p\geq 3$. So let us assume that $p=2$.
$${\rm H}^q(\fom{2})= \liminv{n}{\rm H}^q(K_{\mid nD})$$
and since $K_{\mid D}=-D_{\mid D}$ is non-torsion we get the result.
The same reasoning also gives,
$$E_1^{0,0}={\rm H}^0(\fo)=\C\quad E_1^{0,1}={\rm H}^1(\fo)=\C.$$

Thus we see that $E_2^{p,q}=E_{\infty}^{p,q}$ for all $p,q$.
Since $E_1^{2,1}=0$ we have an exact sequence,
$$0\to E_2^{0,1}\to E_1^{0,1}\stackrel{d}{\to} E_1^{1,1}\to E_2^{1,1}\to 0.$$
Since ${\rm H}^0(\fo)=\C$ one sees that $E_2^{1,0}=E_1^{1,0}$.
Since ${\rm H}^1(D)=0$, one gets $E_2^{1,0}=0=E_2^{0,1.}$ Thus the map

\begin{equation}
d: {\rm H}^1(\fo)\to{\rm H}^1(\fom{1})\label{equat}
\end{equation}
is injective and ${\rm H}^0(\fom{1})=0$.
Since $E_2^{0,2}=E_2^{2,0}=0$, one sees that ${\rm
H}^2(D,\C)=E_2^{1,1}$ and hence, ${\rm H}^1(\fom{1})=\C^{10}$.

The injectivity of $d$ in \ref{equat} above shows that
$\C={\rm H}^1({\cal O}_{(n+1)D}) \stackrel{d}{\longrightarrow}
{\rm H}^1(\Omega^1_{\mid nD})$ is
injective for large $n$. Since ${\rm H}^1({\cal O}_{(n+1)D})$
parametrises line bundles of degree zero on $(n+1) D$,
we see that  any non-trivial
line bundle of degree zero on $(n+1)D$ has non-zero class in
${\rm H}^1(\Omega^1_{\mid nD})$. Thus for large $n$ the class of $D$
in ${\rm H}^1(\Omega^1_X)$ goes to a non-zero element in ${\rm
H}^1(\Omega^1_{\mid nD})$. I claim that this implies the map ${\rm
H}^1(\Omega^1_X)\to {\rm H}^1(\Omega^1_{\mid nD})$ is injective.
Since $X$ is rational we may find a basis for ${\rm
H}^1(\Omega^1_X)$ by the classes of elements of ${\rm Pic}\ X$.
So if $E_i$ form the components of $D$ and $H$ is a general
hyperplane of $X$, then their classes form a basis for ${\rm
H}^1(\Omega^1_X)$. Assume that an element $L=\sum x_iE_i+xH$ goes
to zero in ${\rm H}^1(\Omega^1_{\mid nD})$. One easily sees that
the degree of this element should be zero restricted to any component of
$D$. Thus $L\cdot D=0$ which implies
$x=0$. Now since $L\cdot E_i=0$, one finds that $L=aD$ for
some complex number $a$. But we know that no non-zero multiple of $D$
goes to zero and so $a=0$.

Thus ${\rm H}^0(\Omega^1_{\mid nD})\onto {\rm H}^1(\Omega(-nD))$ for
all large $n$. But ${\rm H}^0(\fom{1})=0$ implies
$\liminv{n}{\rm H}^1(\Omega(-nD))=0$ and by duality,
$\limdir{n}{\rm H}^1(\Omega(nD))=0$. Thus by lemma [\ref{three}], ${\rm
H}^1(\Omega^1_Y)=0$.
\end{proof}

We finally exhibit such surfaces as in the Theorem. The
existence of such surfaces with no condition on $D_{\mid D}$ is
standard. We just follow a similar recipe for our case too.

First we construct such surfaces with $D+K=0$. In the case of
$\widetilde{E}_8$ we take a smooth elliptic curve $C$ and $L$ a flex,
in $\P^2$. Taking the pencil of curves given by $C$ and $3L$ and
resolving its base points one sees that the member of the pencil
containing the proper transform $D$ of $3L$ is of $\widetilde{E}_8$
type and $D+K=0$. In the case of $\widetilde{D}_8$ we start with a
smooth quadric $Q$ and a tangent $L$ to it in $\P^2$. By Bertini, one
sees that for a general linear form $M$ the curve $C=L^3+MQ$ is smooth
elliptic. Take the linear pencil given by $C$ and $LQ$ and resolve
its base points. We end up with a surface with a fiber $D$ of the pencil
of the required type and $D+K=0$.

Unfortunately in both the above cases $D_{\mid D}$ is trivial. Next
we modify these surfaces. Since $D+K=0$, there is an exceptional
curve meeting $D$ exactly once. Blow this down to obtain a surface
with $D'+K$ still zero where $D'$ is the image of $D$ but $K^2=1$. It
is clear that $h^0(D')=2$ and let $P$ be the base point of this
linear system. If we blow up $P$ we will get back to the surface we
started with. Let $E$ be the irreducible component of $D$ containing
$P$. Then $E^2=-1$ and it occurs with multiplicity one in $D$. Let
$Q\neq P$ be any point of $E$ which is smooth on $D$. (This means
that we have to avoid the unique pont of intersection of $E$ with the
rest of the components of $D$). Blow up $Q$ and call the new divisor
$D''$. I claim that $D''$ has all the required properties. Notice that
by choice of $Q$, $h^0(D'')=1$.

Clearly $D''$ has the appropriate configuration and $D''+K=0$.
The only
thing we need to verify is that $D''_{\mid D''}$ is non-torsion. Notice
that ${\rm Pic} D''={\bf G}_a$ and hence $D''_{\mid D''}$ is torsion
is equivalent to it being trivial. Also $D''_{\mid D''}$ is trivial
is equivalent to $h^0(D'')=2$ which is not the case. This completes the
construction.

\begin{rmk}
In the construction above we saw that for each choice of $Q$, we got
a surface with the required properties. I do not know whether these
surfsces are isomorphic. In other words does there exists a family
parametrised by ${\bf G}_a$ of such surfaces?
\end{rmk}

\end{document}